\documentclass[a4paper,12pt]{article}

\usepackage[austrian,american]{babel}

\usepackage[intlimits,leqno]{amsmath}
\usepackage{amssymb} 
\usepackage{amsthm}
\usepackage{graphicx,color}
\usepackage{makeidx}
\usepackage{epsfig}
\usepackage{latexsym}
\usepackage{pstricks,pst-node,amsfonts}

\theoremstyle{plain}

\newcommand\R{{\mathbb{R}}}
\newcommand\C{{\mathbb{C}}}

\newcommand{\F}{\mathcal{F}}
\newcommand{\I}{\mathcal{I}}
 
\renewcommand{\i}{{\rm i}}

\newcommand{\x}{\mathbf{x}}
\renewcommand{\k}{\mathbf{k}}

\renewcommand{\d}{{\rm d}}

\begin{document}

\title{Time reversal for photoacoustic tomography based on the wave equation of Nachman, Smith and Waag}

\author{ Richard Kowar\\
Department of Mathematics, University of Innsbruck, \\
Technikerstrasse 21a, A-6020, Innsbruck, Austria
}

\maketitle

\begin{abstract}
The goal of \emph{photoacoustic tomography} (PAT) is to estimate an \emph{initial pressure function} $\varphi$ from 
pressure data measured at a boundary surrounding the object of interest. This paper is concerned with a time reversal 
method for PAT that is based on the dissipative wave equation of Nachman, Smith and Waag~\cite{NaSmWa90}. This equation 
has the advantage that it is more accurate than the \emph{thermo-viscous} wave equation. For simplicity, we 
focus on the case of one \emph{relaxation process}. 
We derive an exact formula for the \emph{time reversal image} $\I$, which depends on the \emph{relaxation time} $\tau_1$ 
and the \emph{compressibility} $\kappa_1$ of the dissipative medium, and show $\I(\tau_1,\kappa_1)\to\varphi$ for 
$\kappa_1\to 0$. This implies that $\I=\varphi$ holds in the dissipation-free case and that $\I$ is similar to $\varphi$ 
for sufficiently small compressibility $\kappa_1$. Moreover, we show for tissue similar to water that the 
\emph{small wave number approximation} $\I_0$ of the time reversal image satisfies $\I_0 = \eta_0 *_\x \varphi$ 
with $\hat \eta_0(|\k|)\approx const.$ for $|\k|<< \frac{1}{c_0\,\tau_1}$.  
For such tissue, our theoretical analysis and numerical simulations show that the time reversal image $\I$ is very similar 
to the initial pressure function $\varphi$ and that a resolution of $\sigma\approx 0.036\cdot mm$ is feasible 
(in the noise-free case). 
\end{abstract}

\section{Introduction}\label{sec-intro}

The enhancement of \emph{photoacoustic tomography} (PAT) by taking dissipation into account is currently a very active subfield 
of PAT with various contributions from engineers, mathematicians and physicists. 
For basic facts on PAT in the absence and presence of dissipation, we refer for example  
to~\cite{BurMatHalPal07,FinPatRak04,HaScBuPa04,KucKun08,ScGrLeGrHa09,XuWan05,XuWan05b,XuWanAmbKuc03} 
and~\cite{AmBrGaWa11,BurGruHalNusPal07,HriKucNgu08,KaSc13,KoSc12,RivZhaAna06,TrZhCo10}, respectively. 
There are several strategies for solving this class of problem. Firstly, \emph{regularization methods}  
focus on a rigorous mathematical analysis and numerical solution of the inverse problem bearing its degree of 
ill-posedness in mind. Due to its mathematical character, this aspect of PAT is usually carried out by mathematicians. 
This work does not deal with this aspect of PAT. 
Secondly, exact and approximate reconstruction formulas and time reversal methods are developed and investigated by engineers, 
mathematicians and physicists. It is beneficial that these methods furnish physical insights of the considered inverse problem. 
This paper focuses on physical and arithmetical aspects of time reversal of dissipative pressure waves satisfying the 
\emph{wave equation of Nachman, Smith and Waag}  (cf.~\cite{NaSmWa90,KoSc12}) endowed with the respective source term of PAT. 
This equation has the advantage that it is more accurate than the \emph{thermo-viscous} wave equation and permits the modeling 
of several relaxation processes. 
Recently, a \emph{small frequency approximation} of a time reversal functional for PAT based on the 
\emph{thermo-viscous wave equation} were proposed and investigated in~\cite{AmBrGaWa11}. Subsequently, motivated by the 
work~\cite{AmBrGaWa11}, a \emph{small frequency approximation} of a time reversal functional that is based on the 
\emph{wave equation of Nachman, Smith and Waag} was derived in~\cite{KaSc13}. 
The goal of the following paper is to derive and investigate an \emph{exact time reversal formula} as well as its  
\emph{small wave number approximation}. 
In contrast to~\cite{AmBrGaWa11,KaSc13}, we carry out our investigation in the \emph{wave vector-time domain}, which  
permits a more detailed analysis. For simplicity, we concentrate on the case of one relaxation process.

\subsection*{The wave equation model}

In order to explain the results of this paper in more details, we start with the wave equation model which reads for one relaxation 
process as follows (cf. page 113 in~\cite{KoSc12})
\begin{equation}\label{weqNSW}
     \left(\mbox{Id} + \tau_1\frac{\partial}{\partial t}\right)\Delta p 
     - \frac{1}{c_0^2}\left(\mbox{Id} + \tau_0\frac{\partial}{\partial t}\right)\frac{\partial^2 p}{\partial t^2} 
     = - \frac{\varphi}{c_0^2}\left(\mbox{Id} + \tau_1\frac{\partial }{\partial t}\right)
            \delta'(t) \,,
\end{equation} 
where $\varphi$ and $\tau_1>0$ correspond to the \emph{initial pressure function} and the \emph{relaxation time}, respectively. 
If $c_\infty$, $\rho$ and $\kappa_1$ denote the \emph{speed of sound}, the \emph{density} and the \emph{compressibility} of the 
medium, respectively, then 
\begin{equation}\label{deftau0}
   \tau_0 = (1-c_0^2\,\rho\,\kappa_1)\,\tau_1 > 0  \qquad\mbox{and}\qquad c_\infty = \sqrt{\frac{\tau_1}{\tau_0}}\,c_0\,. 
\end{equation}
As shown in~\cite{NaSmWa90,KoSc12}, $\tau_0\in (0,\tau_1)$ is required to guarantee a finite wave front speed $c_F$. 
Then $c_F \leq c_\infty$ holds. 

The direct problem associated to this equation consists in calculating the pressure function $p$ for given initial pressure 
data $\varphi$. Now to the inverse PAT problem.

\subsection*{The time reversal method}

The PAT problem considered in this paper consists in the estimation of the initial pressure function $\varphi$ from pressure 
data $p$ measured at a boundary $\partial \Omega$ that satisfies wave equation~(\ref{weqNSW}). Here it is tacitly assumed that 
$\varphi$ vanishes outside $\Omega$ and $T>0$ is a sufficiently large time period such that every $\x\in\Omega$ satisfies 
$|\x|< c_0\,T\,$. 
The algorithm of our time reversal method can be described as follows: First calculate the function 
$$
       \phi_T = p|_{t=T}   \qquad \mbox{on a sufficiently large domain $\Omega_0$ containing $\Omega$}
$$
from the boundary data on $\partial \Omega$ and the pressure data inside $\Omega$ at time $T$. The latter data are negligible 
for sufficiently large $T$. Secondly, solve the time reversed wave equation
\begin{equation}\label{trweqNSW}
\begin{aligned}
     \left(\mbox{Id} - \tau_1\frac{\partial}{\partial t}\right)\Delta\,q 
     - \frac{1}{c_0^2}\left(\mbox{Id} - \tau_0\frac{\partial}{\partial t}\right)\frac{\partial^2 q}{\partial t^2} 
     = \frac{\phi_T}{c_0^2}\left(\mbox{Id} - \tau_1\frac{\partial }{\partial t}\right)\delta'(t) 
\end{aligned}
\end{equation} 
on $\Omega_0\times [0,T]$ with the help of the calculated function $\phi_T$. Finally, define the 
\emph{time reversal imaging functional} by\footnote{The factor $2$ occurs in the functional, because due to the source term 
in~(\ref{trweqNSW}) only half of the energy of the original pressure wave propagates inside $\Omega$. The other half propagates 
outside $\Omega$.}
\begin{equation}\label{defimagf}
     \phi_T \mapsto \I[\phi_T] := 2\,q|_{t=T} \,.
\end{equation} 
In this paper we derive from the exact representation of $\I$ that
\begin{equation}\label{propF1}
   \lim_{\kappa_1\to 0} \I
           = \varphi 
           = \I|_{\tau_1=0}  \qquad\mbox{, in particular}\qquad
   \I \approx  \varphi  \,  
     \quad \mbox{for small $\kappa_1$}\,.
\end{equation}
This result justifies our time reversal method. We emphasize that 
$\I=\varphi$ cannot hold for dissipative media ($\kappa_1\not=0$). 
Moreover, we show for tissue similar to water that the \emph{small wave number approximation} $\I_0$ of the imaging functional satisfies 
\begin{equation*}
\begin{aligned}  
    \I_0 = \eta *_\x \varphi
       \qquad\mbox{with}\qquad \hat\eta\approx\hat\eta(|\mathbf{0}|) \qquad \mbox{for}\qquad  |\k|<<k_c:=\frac{2}{c_0\,\tau_1}\,,
\end{aligned}           
\end{equation*}
where $\hat\eta(0)\to 1$ for $\kappa_1\to 0$. 
Here $\hat \eta$ denotes the \emph{Fourier transform} of $\eta$ with respect to $\x\in\R^3$ and 
$*_\x$ denotes the \emph{space convolution}. 
In addition, our analysis and simulations indicate that a resolution of $\sigma\approx 0.036\cdot mm$ is 
feasible (in the noise-free case). 

We note that the theoretical results of this paper are valied for any space dimension $\geq 1$. \\

This paper is organized as follows. In Section~\ref{sec-opeq} we derive an exact representation of $\I$ and discuss the 
limit $\kappa\to 0$ for vanishing compressibility. 
The small wave number representation of $\I$ is derived and discussed in the subsequent section. 
For a better understanding of the results, we have visualized many of our results (like eigenvalues, coefficients and so on) 
with MATLAB. 
Moreover, for the convenience of the reader, we have listed the solution of the dissipative wave equation and its time 
reversal in the appendix. 
Finally, the paper is concluded with the section Conclusions.

\section{Time reversal image $\I$ and its properties}
\label{sec-opeq}

First we derive an exact representation of the time reversal image $\I$ defined by~(\ref{defimagf}). 
Throughout this paper we denote by $k:=|\k|$ the wave number corresponding to the wave vector $\k\in\R^3$.  
Let 
\begin{equation*}
\begin{aligned}  
        \hat\phi_T(\k) :&= \hat p(\k,T) = \hat\varphi(\k)\,\sum_{j=0}^2  A_j(k)\,\lambda_j(k)\,e^{-\lambda_j(k)\,T}  
                   \qquad\mbox{and}\\
        \hat q(\k,T)    &= - \hat\phi_T(\k)\,\sum_{l=0}^2  A_l(k)\,\lambda_l(k)\,e^{\lambda_l(k)\,T}
\end{aligned}           
\end{equation*}
be the solutions of the wave equation~(\ref{weqNSW}) and its time reverse~(\ref{trweqNSW}) at time $t=T$, respectively, 
derived in the appendix. The representations for $A_j$ and $\lambda_j$ for $j=0,\,1,\,2$ are listed in the appendix and 
visualized in Fig.~\ref{fig:A} and Fig.~\ref{fig:01}. In this paper we focus on tissue similar to water  
for which $\lambda_0=\lambda_0(k)$ is real-valued and
$$
             \lambda_1 = \mu + \i\,\vartheta  \qquad\mbox{and}\qquad \lambda_2 = \mu - \i\,\vartheta
$$
with real-valued $\mu=\mu(k)$ and $\vartheta=\vartheta(k)$ listed in the appendix (cf. Fig.~\ref{fig:01}). 
\begin{figure}[!ht]
\begin{center}
\includegraphics[height=5cm,angle=0]{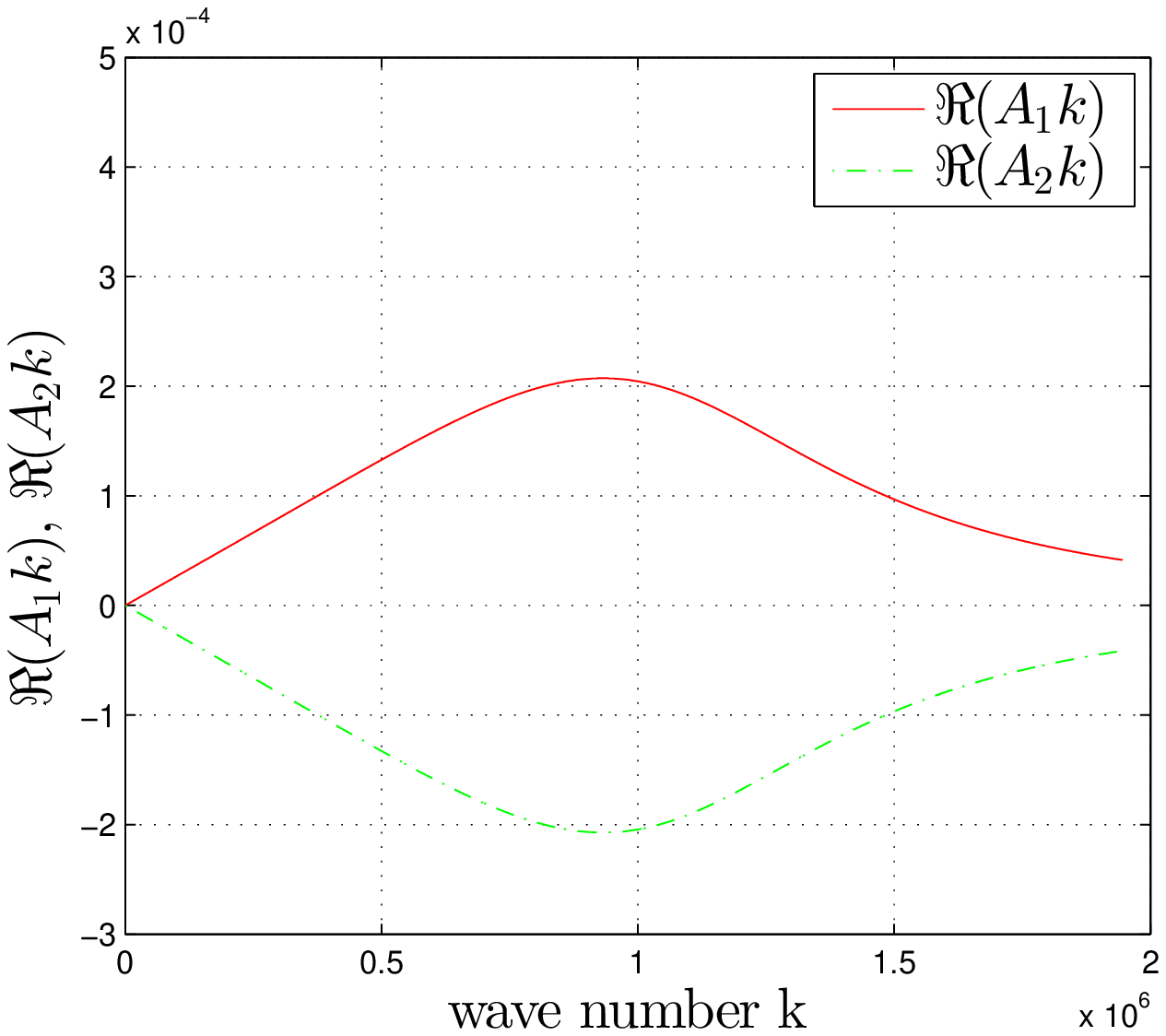}
\includegraphics[height=5cm,angle=0]{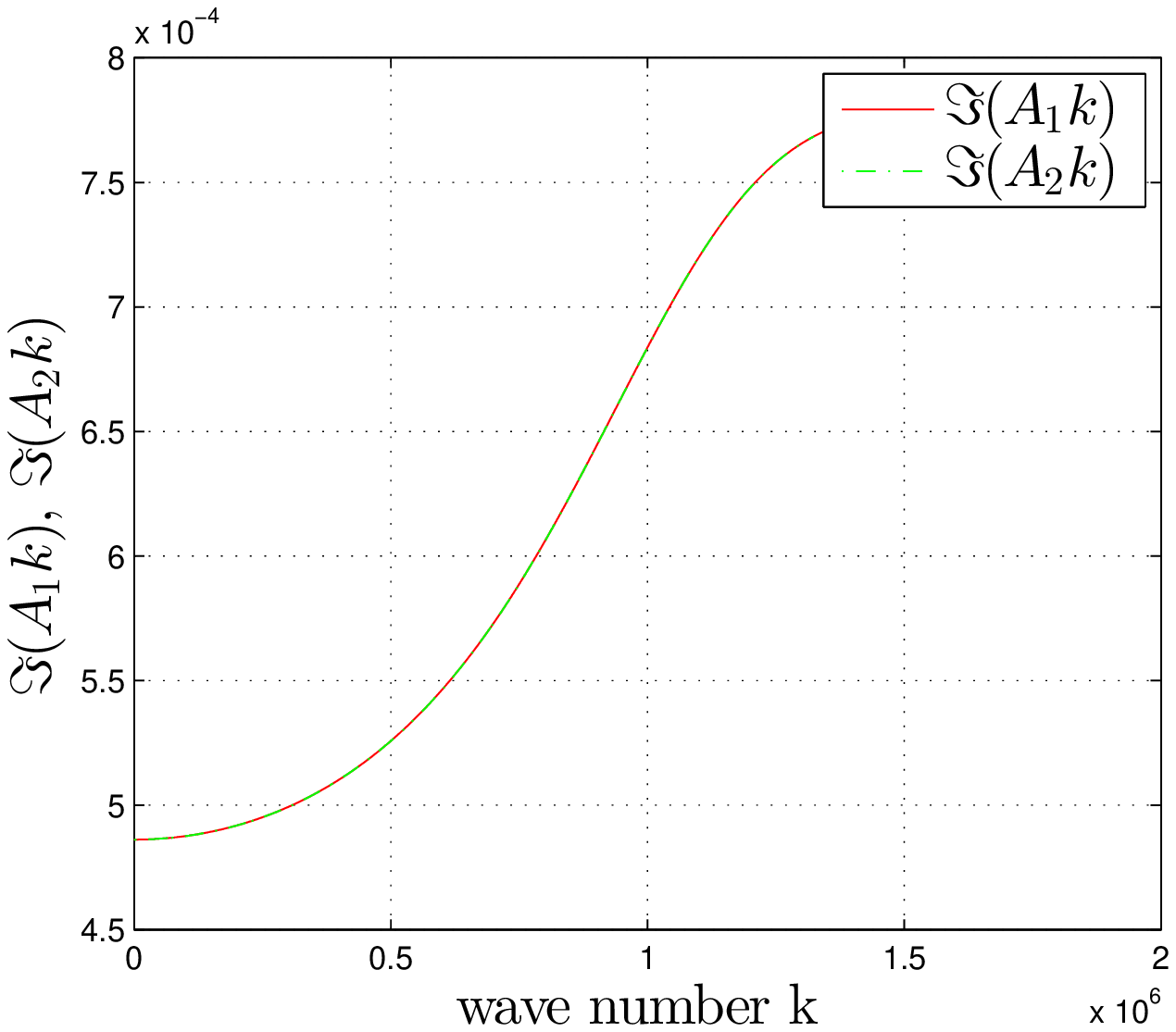}\\
\includegraphics[height=5cm,angle=0]{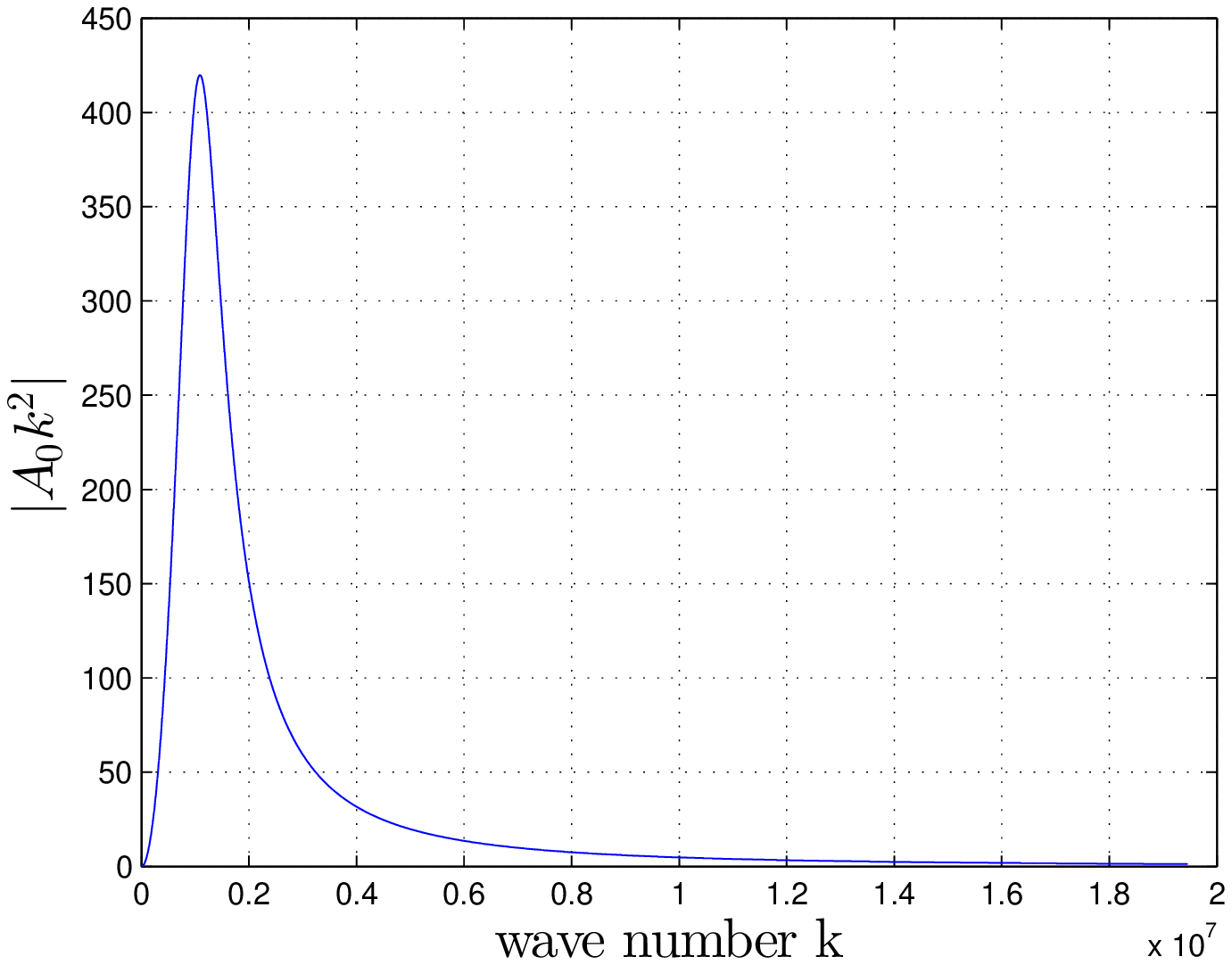}
\includegraphics[height=5cm,angle=0]{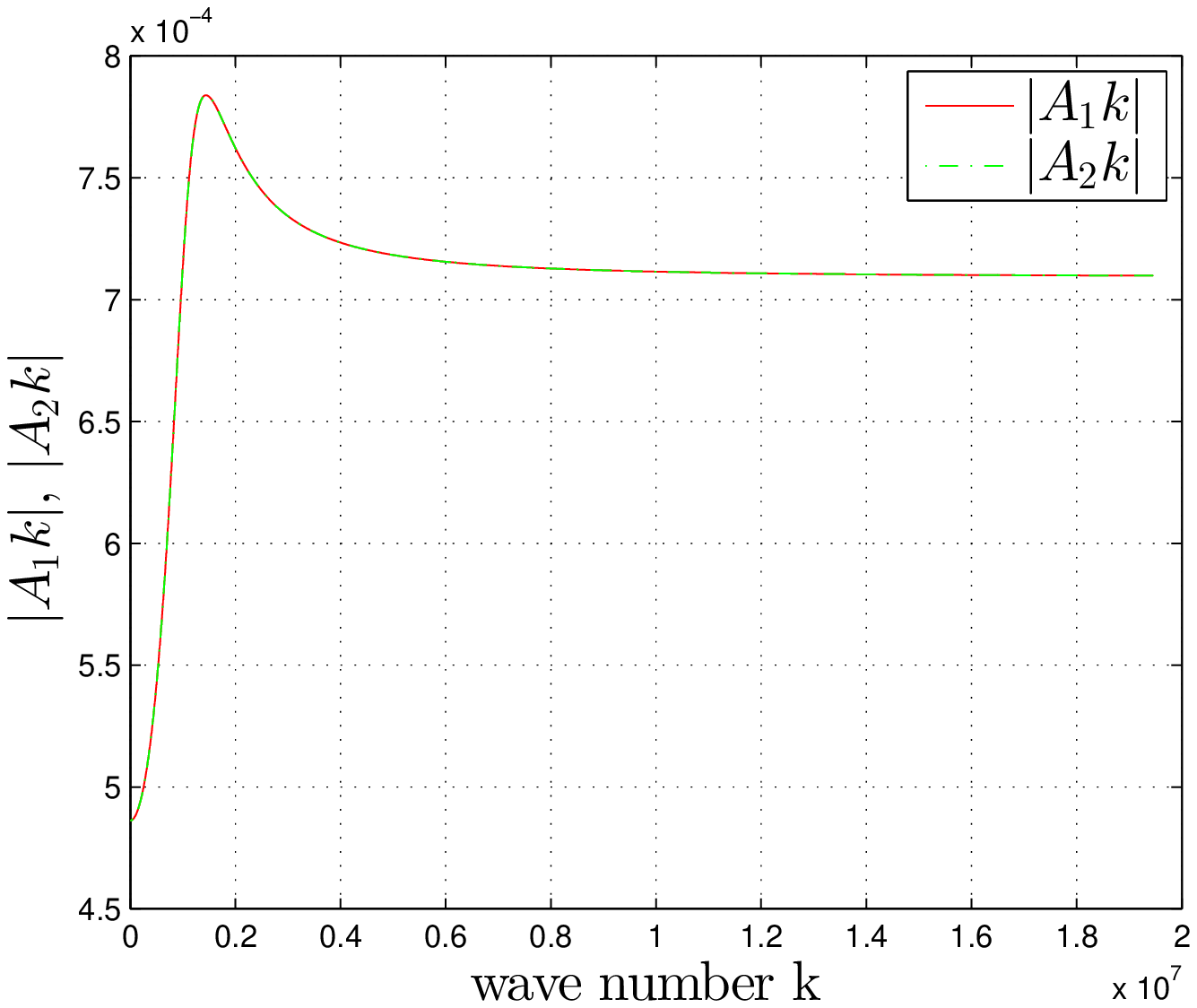}
\end{center}
\caption{The upper row visualizes $\Re(A_j*k)$ and $\Im(A_j*k)$ ($j=1,\,2$) for $k\in [0,k_c]$. We see that $\Re(A_1)=-\Re(A_2)$ 
and $\Im(A_1)=\Im(A_2)$.  
The lower row visualizes $|A_0*k^2|$, $|A_1*k|$ and $|A_2*k|$ for $k\in [0,10\cdot k_c]$. We see that $|A_1|=|A_2|$ and each 
$|A_j|$ has the same growth behavior as in~(\ref{O2}).}
\label{fig:A}
\end{figure}

Employing these representations and $\cosh(z) = \frac{1}{2}\,(e^{z} + e^{-z})$ for $z\in\C$ to~(\ref{defimagf}) results in 
\begin{equation*} 
       \hat \I 
          =   \left[2\,\sum_{j=0}^2 A_j^2\,\lambda_j^2
            + 4\,\sum_{j,l=0,\,j<l}^2  A_j\,A_l\,\lambda_j\,\lambda_l\,\cosh((\lambda_j-\lambda_l)\,T)\right]\,\hat\varphi\,.       
\end{equation*}
Because $\lambda_2 = \overline{\lambda_1}$, $A_2 = -\overline{A_1}$ and $\lambda_1-\lambda_2=\i\,2\,\vartheta$ (cf. Appendix), 
it follows with 
\begin{equation}\label{propcosh2}
 \cosh(x+\i\,y) = \cosh(x)\,\cos(y) +\i\,\sinh(x)\,\sin(y)
\end{equation}
and 
\begin{equation}\label{propcosh}
          \cosh(\i\,2\,x) = \cos(2\,x) = 1-2\,\sin^2(x)
\end{equation}
that 
$$
        \hat \I  = (2\,\pi)^{3/2}\,(\hat\zeta_1 - \hat\zeta_2 + \hat\zeta_3) \,\hat\varphi
$$ 
with \emph{real-valued} functions 
\begin{equation}\label{zeta12}
\begin{aligned}  
   \hat \zeta_1  
       = \frac{2\,\sum_{j=0}^2 A_j^2\,\lambda_j^2 + 4\,|A_1|^2\,|\lambda_1|^2}{(2\,\pi)^{3/2}} \,,\qquad\quad
   \hat \zeta_2
       = \frac{8\,|A_1|^2\,|\lambda_1|^2\,\sin^2(\vartheta\,T)}{(2\,\pi)^{3/2}}  
\end{aligned}           
\end{equation}
and 
\begin{equation}\label{zeta3}
\begin{aligned}  
    \hat \zeta_3
       &= \frac{8\,A_0\,\lambda_0}{(2\,\pi)^{3/2}}\,\left[
           \Re(A_1\,\lambda_1)\,\cosh(\Re(\lambda_0-\lambda_1)\,T)\,\cos(\Im(\lambda_0-\lambda_1)\,T) \right.\\
          &\qquad\qquad \quad \left. - \Im(A_1\,\lambda_1)\,\sinh(\Re(\lambda_0-\lambda_1)\,T)\,\sin(\Im(\lambda_0-\lambda_1)\,T)
       \right]  \,.
\end{aligned}           
\end{equation}
For the last identity we used 
$$
       A_1\,\lambda_1+A_2\,\lambda_2 = 2\,\Re(A_1\,\lambda_1) \qquad\mbox{and}\qquad 
       A_1\,\lambda_1-A_2\,\lambda_2 = \i\,2\,\Im(A_1\,\lambda_1)\,.
$$
Applying the inverse Fourier transform and the convolution theorem (cf. Appendix) to the image $\I$ yields
\begin{equation}\label{repI}
\begin{aligned}  
    \I  = (\zeta_1 - \zeta_2 + \zeta_3) *_\x \varphi \,,
\end{aligned}           
\end{equation}
where $*_\x$ denotes the convolution with respect to the space variable $\x$. 
\begin{figure}[!ht]
\begin{center}
\includegraphics[height=5cm,angle=0]{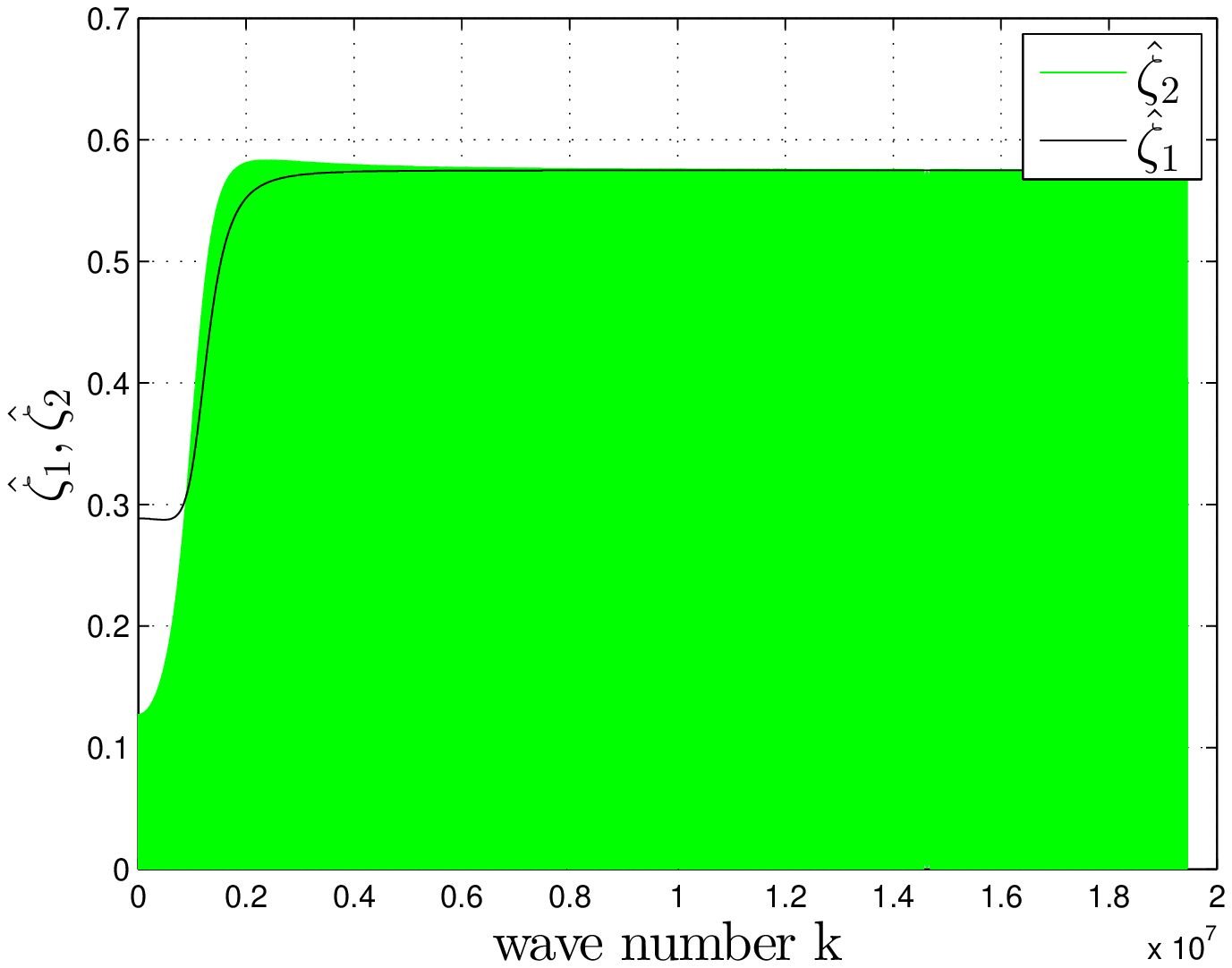}
\includegraphics[height=5cm,angle=0]{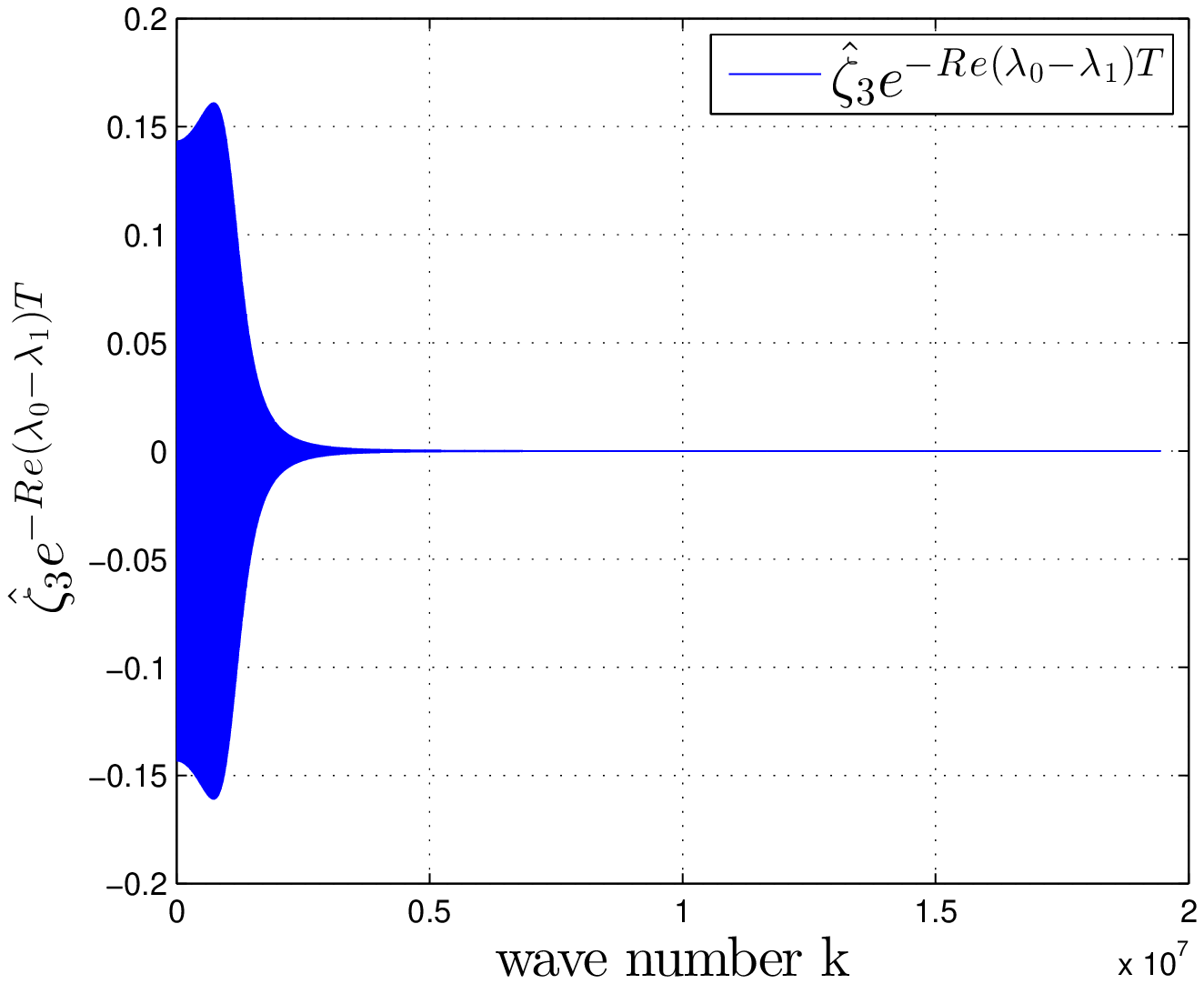}
\end{center}
\caption{Visualization of $\hat\zeta_j$ for $j=0,\,1,\,2$. 
We see that $\hat\zeta_2$ and $\hat\zeta_3$ are highly oscillating and that each $\hat\zeta_j$ has the growth behavior 
as in~(\ref{O3}). Although $\Re(\lambda_0-\lambda_1)$ is bounded, it is of order $10^9$, i.e. MATLAB cannot calculate 
$\cosh(\Re(\lambda_0-\lambda_1)\,T)$ and thus $\hat\zeta_3\cdot e^{-\Re(\lambda_0-\lambda_1)\,T}$ was plotted  instead of 
$\hat\zeta_3$.}
\label{fig:zeta}
\end{figure}
From the representations of $\lambda_j$ and $A_j$ for $j=0,\,1,\,2$ in the appendix, it follows that (cf. Fig.~\ref{fig:01})
\begin{equation}\label{O1}
     \lambda_0(k),\,\mu(k) = \mathcal{O}(1) , \qquad
     \vartheta(k) = \mathcal{O}(k) 
     \qquad\mbox{for}\qquad 
     k\to\infty
\end{equation}
and (cf. Fig.~\ref{fig:A})
\begin{equation}\label{O2}
     |A_0(k)| = \mathcal{O}\left(\frac{1}{k^2}\right)\,,\quad  |A_1(k)|,\,|A_2(k)| = \mathcal{O}\left(\frac{1}{k}\right) 
      \quad\mbox{for}\quad 
     k\to\infty\,.
\end{equation}
In this paper we focus on tissue similar to water for which $\vartheta$ is real-valued,\footnote{We stress this fact/advantage, 
because if the thermo-viscous wave equation is used, then $\vartheta$ is complex-valued for large $k$ and thus $\sin^2(\vartheta\,T)$ 
is exponentially increasing, which implies a stronger restriction on $\varphi$.} 
i.e. $\sin^2(\vartheta\,T)$ 
is bounded and thus $\zeta_2$ is bounded, too. 
From~(\ref{zeta12}),~(\ref{zeta3}) with the boundedness of $\Re(\lambda_0-\lambda_1)$,~(\ref{O1}) and~(\ref{O2}), it follows 
for each $j\in\{0,\,1,\,2\}$ that (cf. Fig.~\ref{fig:zeta})
\begin{equation}\label{O3} 
      |\hat\zeta_j(k)| = \mathcal{O}(1)
      \quad\mbox{for}\quad    k\to\infty \,.          
\end{equation}

\subsection*{The time reversal condition for $\varphi$}

The previous asymptotic relations show that time reversal can be applied successfully if the initial pressure 
function $\varphi$ is a \emph{quadratic integrable} function, because 
$$
   \int |\hat\I(\k)|^2 \,\d \k 
         \leq C\, \int |\hat\varphi(\k)|^2 \,\d \k 
         = C\, \int |\varphi(\x)|^2 \,\d \x
$$
holds with 
$$
       C:=(2\,\pi)^{3/2}\,\max_{k\in [0,\infty]}(|\zeta_1(k)|^2,\,|\zeta_2(k)|^2,\,|\zeta_3(k)|^2)<\infty \,,
$$
where $k=|\k|$. That is to say, the time reversal image $\I$ exists and is \emph{quadratic integrable} if $\varphi$ is 
\emph{quadratic integrable}.

\subsection*{Numerical problem in calculating $\F^{-1}\{\hat\zeta_3\}$}

We would like to make a short remark about the numerical problem of calculating $\zeta_3$ as the inverse Fourier transform of 
$\hat\zeta_3$.  
Although $\Re(\lambda_0-\lambda_1)$ is bounded, it is of order $10^9$ in the small wave number region and thus the term 
$\cosh(\Re(\lambda_0-\lambda_1)\,T)$ in the representation of $\hat\zeta_3$ cannot be calculate by MATLAB. 
In addition, $\hat\zeta_3$ is highly oscillatory.  
All in all, $\F^{-1}\{\hat\zeta_3\}$ cannot be calculated numerically in the form~(\ref{zeta3}). 
However, it will be shown in Section~\ref{sec-small} below that $\zeta_3$ is negligible small under appropriate conditions. 
It is noteworthy that this can be explained due to wave propagation.

\subsection*{The limit case $\kappa_1\to 0$}
\label{sec-lim}

Now we show that $\I$ is a good approximation of $\varphi$ if the compressibility $\kappa_1$ is sufficiently small. More precisely, 
we show 
$$
   \lim_{\kappa_1\to 0} \I = \I|_{\kappa_1=0} =  \varphi \,.
$$

From~(\ref{deftau0}) we conclude $\lim_{\kappa_1\to 0}\frac{\tau_1}{\tau_0}=1$ and $\lim_{\kappa_1\to 0} c_0 = c_\infty$ 
and thus wave equation~(\ref{weqNSW}) becomes for $\kappa_1\to 0$ 
$$
  \left( \mbox{Id} + \tau_1\,\frac{\partial}{\partial t} \right)\left( \,\Delta\,p 
     - \frac{1}{c_0^2}\,\frac{\partial^2 p}{\partial t^2} 
     + \frac{\varphi(\x)}{c_0^2}\,\delta'(t) \right) = 0 \,.
$$
We see that $\kappa_1\to 0$ leads to the dissipation free case. Because $\hat p$ depends continuously on $\tau_0$, $\tau_1$ 
and $c_0$, it remains to show that $\I = \varphi$ for $\kappa_1=0$. For $\kappa_1=0$, the characteristic equation~(\ref{eqlambda}) 
in the appendix simplifies to
$$
    (1 - \tau_1\,\lambda) \,( \lambda^2 + c_0^2\,k^2)  = 0 
            \qquad\mbox{for}\qquad k\in [0,\infty) \,,
$$
which has the solutions 
$$
     \lambda_0 = \frac{1}{\tau_1}\,\qquad\mbox{and}\qquad
     \lambda_{1,2} := \pm \i\,\vartheta  \qquad\mbox{with}\qquad \vartheta:=c_0\,k \,.
$$
Moreover, the expressions for $A_j$ listed in the appendix simplify to  
$$
      A_0=0\,\qquad\mbox{and}\qquad 
      A_1 = -A_2 = - \frac{1}{2\,\lambda_1} \,.  
$$
Employing these results to~(\ref{zeta12}) and~(\ref{zeta3}) yields 
$$
        \zeta_1 = (2\,\pi)^{-3/2}\,2\,,\qquad  
        \zeta_2 = (2\,\pi)^{-3/2}\,2\,\sin^2(c_0\,k\,T)\qquad\mbox{and}\qquad 
        \zeta_3 = 0
$$ 
and therefore 
\begin{equation}\label{limI}
    \hat\I 
       = 2\,\hat\varphi - 2\, \sin^2(c_0\,k\,T)\,\hat\varphi\,.
\end{equation}
In Theorem~3 in~\cite{Ko13}, it is shown that 
\begin{equation}\label{propdelta}
          \F^{-1}\{\sin^2(c_0\,k\,T)\,\hat\varphi\}(\x) = \frac{\varphi(\x)}{2}  \qquad\mbox{for}\qquad 
          \x\in \Omega
\end{equation}
if $T$ is sufficiently large and consequently the claim $\I \to \varphi$ for $\kappa_1\to 0$ follows. 
Here $\F^{-1}$ denotes the inverse Fourier transform.

\subsection*{An interpretation of~(\ref{propdelta})}

It is best to discuss~(\ref{propdelta}) via the equivalent relation 
$$
   -c_0^2\,\Delta\,[ G_0 *_\x G_0](\x,T) = \frac{\delta(\x)}{2}  \qquad\mbox{on $\Omega$ for sufficiently large $T$. }
$$
Here $G_0$ is defined by 
\begin{equation}\label{defG0}
         \hat G_0(\k,t) := (2\,\pi)^{-3/2}\,\frac{\sin(c_0\,k\,t)}{c_0\,k}\,H(t)\, \qquad (k=|\k|,\,\k\in\R^3,\,t\in\R)
\end{equation}
denotes the solution of the \emph{standard wave equation} with source term $f = \delta(\x)\,\delta(t)$ which 
vanishes at time $t=T$ everywhere except at $|\x|=c_0\,T$. Hence $[G_0 *_\x G_0]_{t=T}$ corresponds to a wave at time $t=T$ 
initiated by the source term $f = G_0|_{t=T}\,\delta(t)$ which vanishes everywhere except at those $\x\in\R^3$ satisfying 
$$
        |\x|=0   \qquad\mbox{or}\qquad  |\x|=c_0\,T\,.
$$
In words, one half of the wave initiated at the sphere $|\x|=c_0\,T$ propagates into the sphere and arrives in $\x=\mathbf{0}$ 
at time $T$ and the other half propagates outward the sphere and arrives in $|\x|=2\,c_0\,T$ at time $T$.  
If we focus only on space points in $\Omega$ which do not contain the points $\x$ with $|\x|=2\,c_0\,T$, 
then the above identity follows. If the period $T$ during which the pressure data are acquired is sufficiently large, 
then this assumption is always satiesfied. 

Because~(\ref{propdelta}) for $T/2$ instead of $T$ is required later, we have assumed at the beginning of this paper that $T$ 
is so large that $\Omega\subset\{\x\in\R^3\,|\,|\x|=c_0\,T\}$ holds.

\section{Small wave number approximation}
\label{sec-small}

For theoretical considerations and numerical simulations, it is desirable to have a small wave number approximation of 
the time reversal imaging function $\I$. In the following we derive such an approximation of $\I$ and discuss 
its properties. 

Motivated by the thermo-viscous case, we consider 
$$
        k_c = \frac{2}{c_0\,\tau_1}\,
$$
as the threshold between small and large wave numbers. 
From the representations of $\lambda_j$ for $j=0,\,1,\,2$ in the appendix, it follows for $k << k_c$ that (cf. 
Fig.~\ref{fig:01})
\begin{equation*}
\begin{aligned}
   \lambda_{1,2} \approx \tilde\mu \pm \i\,\tilde\vartheta  \qquad\mbox{with}\qquad 
   \tilde\mu(k) := c_0\,\frac{k^2}{k_c}  \quad\mbox{and}\quad   \tilde\vartheta(k) := c_0\,k\,,
\end{aligned}
\end{equation*}
\begin{equation*}
\begin{aligned}
  \lambda_0(k) 
       \approx 
        \tilde\lambda_0(k) 
        := \frac{1}{\tau_0} -c_0^2\,\tau_1\,k^2  \qquad\mbox{and}\qquad 
      0 < \tilde\mu(k),\,\tilde\vartheta(k)  << \tilde\lambda_0(k) \,.
\end{aligned}
\end{equation*}
\begin{figure}[!hb]
\begin{center}
\includegraphics[height=5.1cm,angle=0]{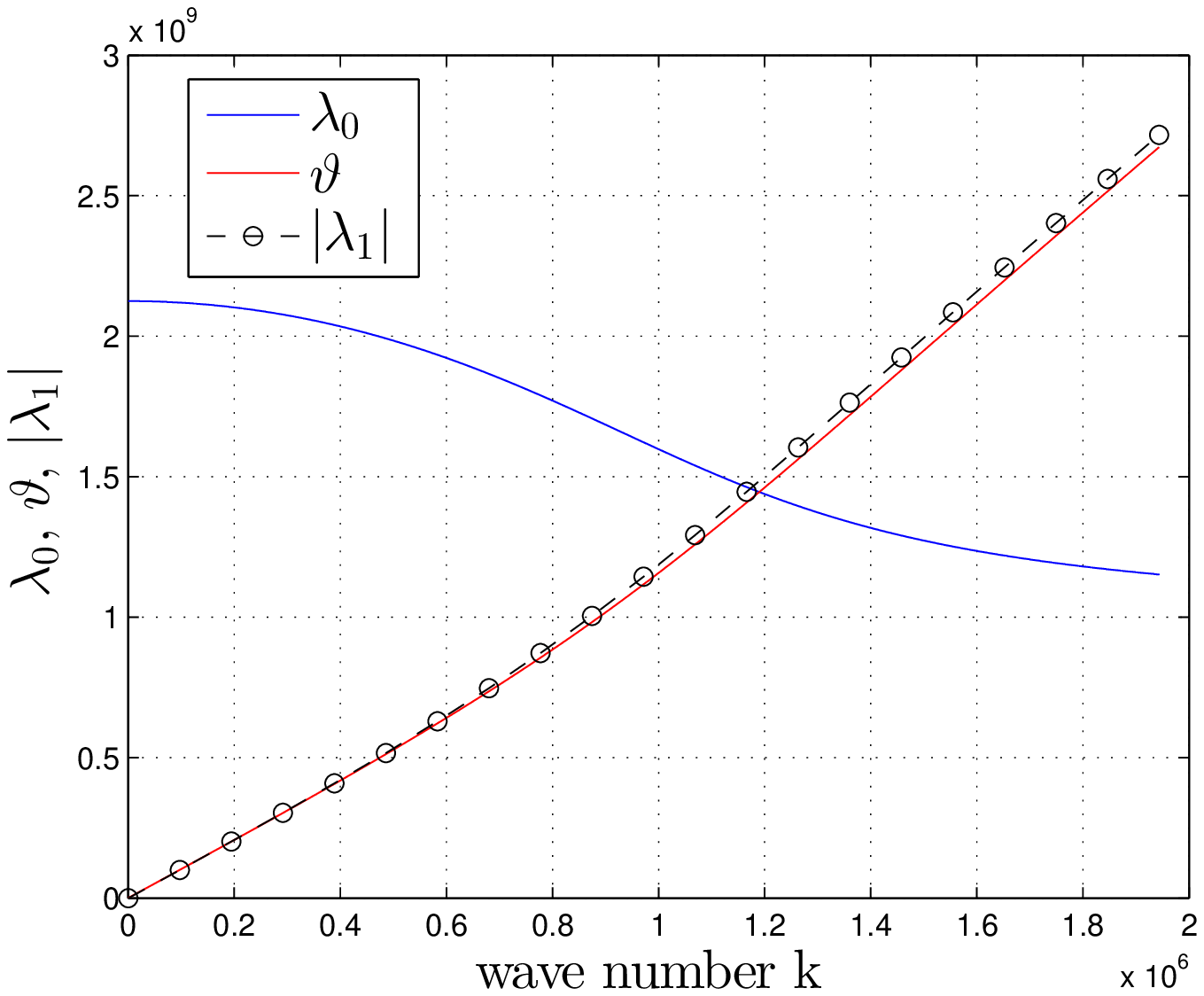}
\includegraphics[height=5.1cm,angle=0]{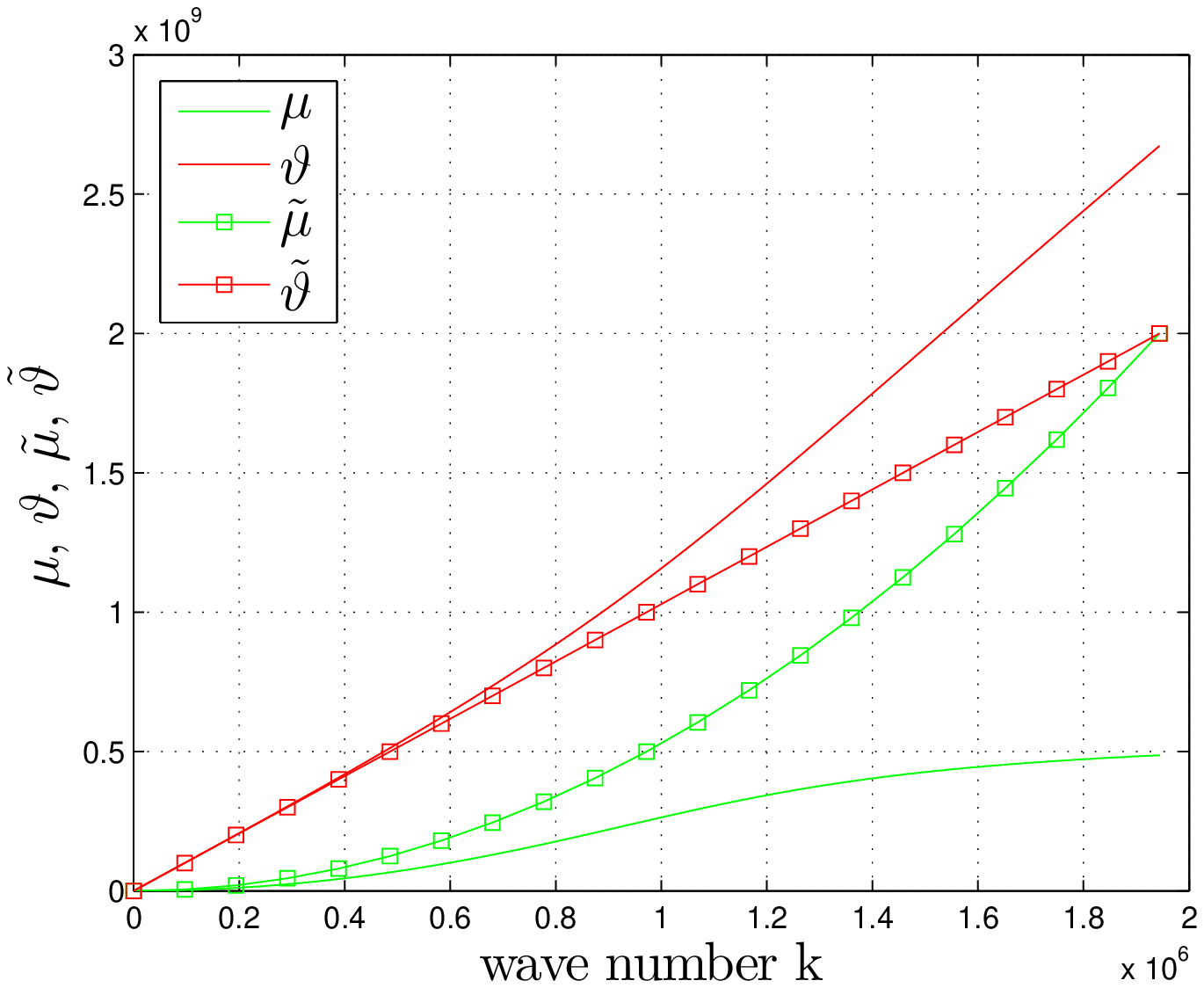}\\
\includegraphics[height=5.1cm,angle=0]{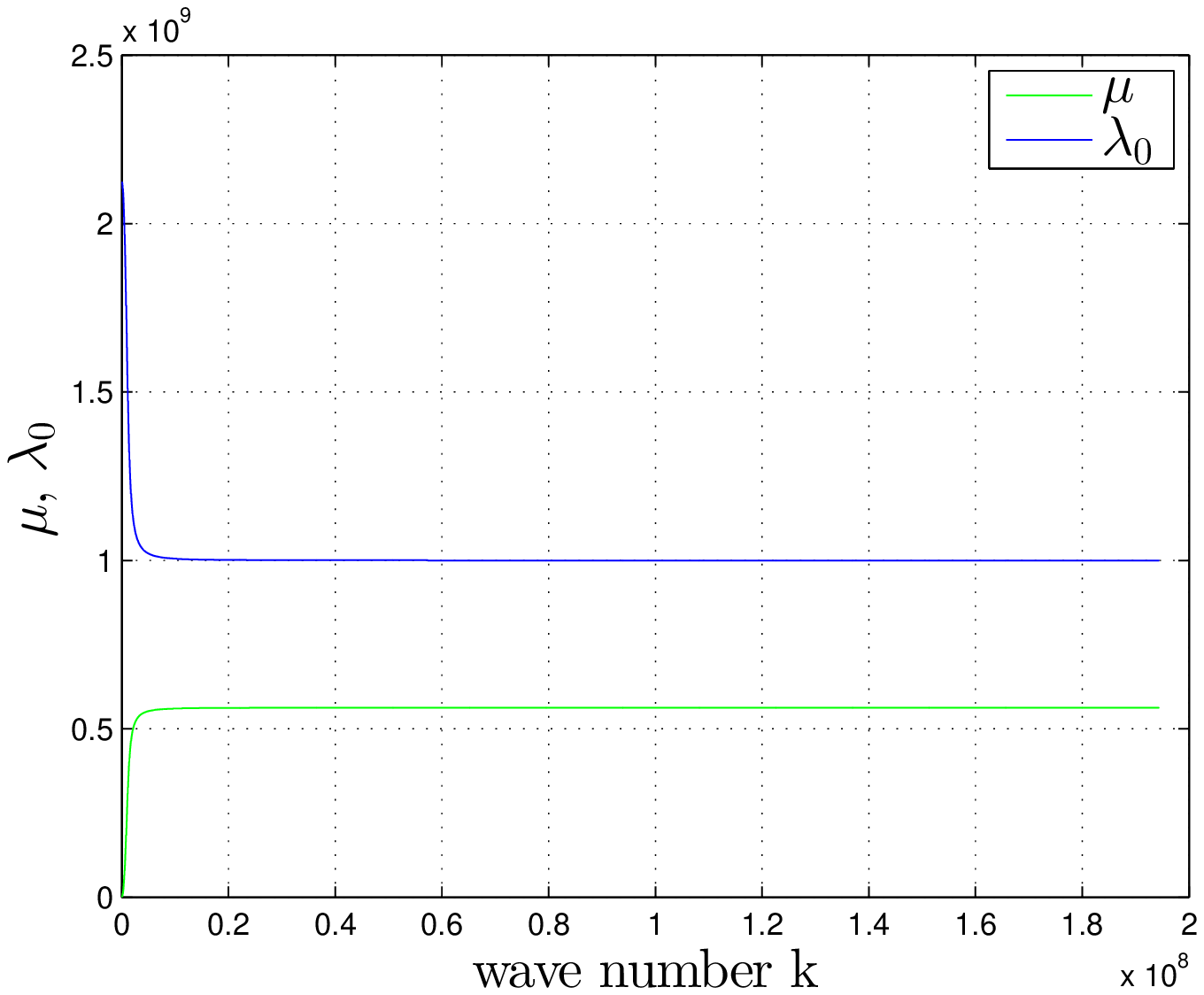}
\includegraphics[height=5.1cm,angle=0]{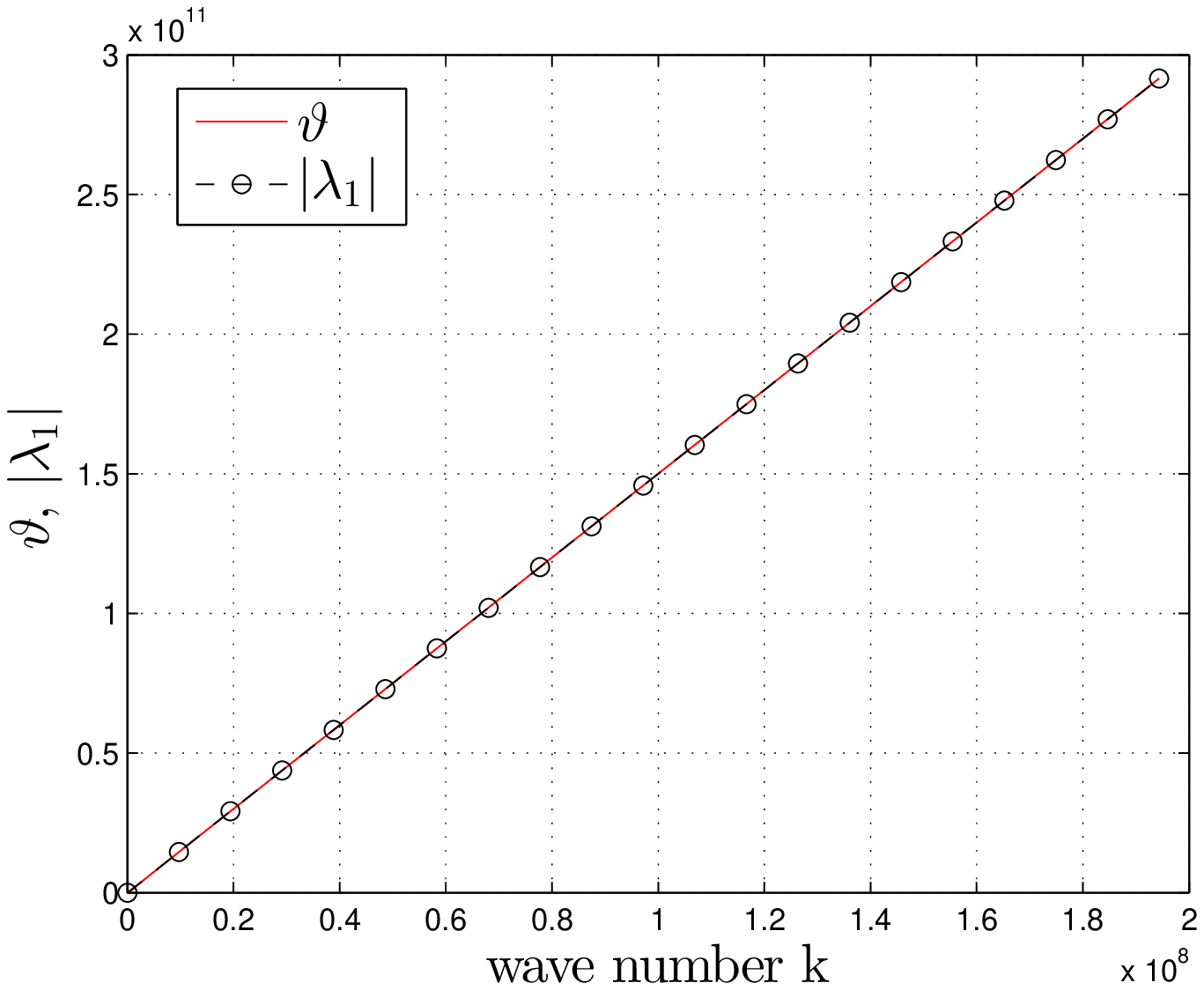}
\end{center}
\caption{The upper row visualizes $\lambda_0$, $\mu$, $\vartheta$ and $|\lambda_1|$ for $k\in [0,k_c]$.  
$\tilde \mu$ and $\tilde\vartheta$ are the small wave number approximations of $\mu$ and $\vartheta$. 
The lower row visualizes $\lambda_0$, $\mu$, $\vartheta$ and $|\lambda_1|$ for $k\in [0,100\cdot k_c]$. We see that 
$\mu\to 5.6250\cdot 10^8\cdot s^{-1}$, $\lambda_0\to 10^9\cdot s^{-1}$ and $\frac{|\lambda_1|}{\vartheta}\to 1$ for 
$k\to\infty$ which agrees with~(\ref{limmu}).}
\label{fig:01}
\end{figure}
This, the convolution theorem (cf. Appendix) and the identity~(\ref{propdelta}) imply 
$$
   \zeta_2 *_\x \varphi 
      \approx \F^{-1}\{8\,|A_1|^2\,|\tilde\lambda_1|^2\,\sin^2(c_0\,k\,T) \,\hat \varphi\}
      =       \F^{-1}\left\{\frac{8\,|A_1|^2\,|\tilde\lambda_1|^2}{(2\,\pi)^{3/2}}\right\} *_\x \frac{\varphi}{2}  \,,
$$
if $T$ is sufficiently large. Hence
\begin{equation}\label{defeta0}
    (\zeta_1+\zeta_2) *_\x \varphi 
         \approx \eta_0 *_\x \varphi  \qquad\mbox{with}\qquad 
     \hat\eta_0 := (2\,\pi)^{-3/2}\,\,2\,\sum_{j=0}^2 A_j^2\,\lambda_j^2\,. 
\end{equation}
Because $\Re(\lambda_0-\lambda_1)\,T$ is very large, we have 
$$
    \cosh(\Re(\lambda_0-\lambda_1)\,T) 
       \approx \sinh(\Re(\lambda_0-\lambda_1)\,T) 
       \approx \frac{1}{2}\,e^{\Re(\lambda_0-\lambda_1)\,T} \,
$$
and thus 
\begin{equation*} 
    \hat \zeta_3 \approx -\hat\eta_1\,\sin(c_0\,k\,T) + \hat\eta_2\,\cos(c_0\,k\,T)    
\end{equation*}
with 
\begin{equation*}
    \hat\eta_1 := \frac{4\,A_0\,\lambda_0}{(2\,\pi)^{3/2}}\,e^{\Re(\lambda_0-\lambda_1)\,T}\,\Im(A_1\,\lambda_1)
        \qquad\mbox{and}  \qquad 
    \hat\eta_2 := \hat\eta_1\,\frac{\Re(A_1\,\lambda_1)}{\Im(A_1\,\lambda_1)}\,.
\end{equation*}
 For the following argumentation, it is required that 
\begin{equation}\label{ass}
        \varphi_\eta (\x) := (\eta_2*_\x \varphi)(\x) \approx 0 \qquad\mbox{for}\qquad   |\x|\geq c_0\,T\,  
\end{equation}
such that~(\ref{propdelta}) holds for $\varphi_\eta$ replacing $\varphi$. 
The crucial point will be that the inverse Fourier transform of the previous sinus and cosinus functions are 
non-dissipative waves, only their coefficients depend on dissipation.
\begin{itemize}
\item [1)] From $\cos(x)=1-2\,\sin^2(x/2)$, $\vartheta\approx c_0\,k$ and identity~(\ref{propdelta}), we get 
           \begin{equation*} 
                 \F^{-1}\{\cos(c_0\,k\,T)\,\eta_2\,\hat\varphi\}
                   \approx \eta_2 *_\x \varphi - \F^{-1}\{2\,\sin^2\left(c_0\,k\,\frac{T}{2}\right)\,\eta_2\,\hat\varphi\} 
                      = 0 
           \end{equation*}
           on $\x\in \Omega$ and for sufficiently large $T$. Hence the term with $\hat\eta_2$ in $\hat\zeta_3$ can be neglected.

\item [2)] Because of 
           $$
               - c_0\,k\,\sin(c_0\,k\,T) = \left.\frac{\partial \cos(c_0\,k\,t)}{\partial t}\right|_{t=T}\,,
           $$
           it follows similarly as in item 1) that
           \begin{equation*} 
                 \F^{-1}\{\sin(c_0\,k\,T)\,\eta_1\,\hat\varphi\} = 0 
           \end{equation*}
           on $\x\in \Omega$ and for sufficiently large $T$. Therefore $\hat\zeta_3$ can be neglected.
           
\end{itemize}
In summary, we get under the assumption~(\ref{ass}) the \emph{small wave number approximation} 
\begin{equation}\label{I1} 
     \I \approx  \I_0 := \eta_0 *_\x \varphi  
          \qquad\mbox{with $\eta_0$ defined as in~(\ref{defeta0}).} 
\end{equation}

\subsection*{Physical interpretation of the term $\zeta_3$}

Let~(\ref{ass}) be satisfied and 
$$
   u 
    := \F^{-1}\{\cos(c_0\,k\,T)\,\eta_2\,\hat\varphi\} 
     = G_0 *_\x \varphi_\eta    
     \qquad\mbox{with $G_0$ as in~(\ref{defG0}).} 
$$
Because $u$ is nothing else but the solution of the \emph{standard wave equation} with source term 
$f=\varphi_\eta\,\delta(t)$, it follows that 
$$
         u(\x,T) = 0  \qquad \mbox{for} \quad \x\in\Omega  \qquad \mbox{,if $T$ is sufficiently large.}
$$
In words, the waves within $\Omega$ propagated out of $\Omega$ and, due to our assumption, no wave propagated into $\Omega$. 
Hence the term $u$ vanishes. Similarly, one shows that $\F^{-1}\{\eta_1\,\sin(c_0\,k\,T)\,\hat\varphi\}$ is the derivative with 
respect to $T$ of a function vanishing inside of $\Omega$ and consequently $\zeta_3$ is negligible on $\Omega$.

\subsection*{Example for tissue similar to water}

The parameter values for tissue similar to water at normal temperatur are (cf.~\cite{KiFrCoSa00})
\begin{equation}\label{parvalues}
  \tau_1\approx 10^{-9}\,s\,,\quad c_\infty\approx 1500\,\frac{m}{s}\,,\quad \rho\approx 10^3\,\frac{kg}{m^3} 
  \quad\mbox{and}\quad
  \kappa_1 \approx 5\,\cdot 10^{-10}\,\frac{m^2}{N}\,
\end{equation}
for which $\tau_0 = \frac{\tau_1}{1+c_\infty^2\,\rho\,\kappa_1} \approx 4.7\cdot 10^{-10}\, s \approx 0.47\cdot\tau_1$ follows. 
As time period, we have chosen the value $T=\frac{4\,L}{c_\infty}$ with $L=0.5\cdot m$. 
We assume that  
\begin{equation*}
        \hat\varphi(\k),\,\hat\eta_1\,\hat\varphi(\k),\,\hat\eta_2\,\hat\varphi(\k) \approx 0 
                      \qquad\mbox{for}\qquad 
        |\k|\leq \frac{k_c}{100}\,
\end{equation*}
such that the small wave number approximation is applicable for $\k\in\R^3$ with $k=|\k|\leq \frac{k_c}{100}$. 
For this range of $\k$, it follows that $\hat\eta_0(k)\approx \hat\eta_0(0)$ with 
$$
         C 
          := (2\,\pi)^{3/2}\,\hat\eta_0(0) 
           =  2\,\left(1-\frac{\tau_1}{\tau_0}\right)^2 + 1 
         \approx 3.5\,.
$$
Hence we arrive at 
$$
    \I_0 
      \approx C\,\varphi 
      \approx 3.5\cdot\varphi  \,.
$$
We note that $C\to 1$ for $\kappa\to 0$, i.e. $C\,\varphi=\varphi$ holds in the dissipation-free case. \\

A simple example satisfying our assumptions is given by
$$
          \varphi(\x) = (4\,\pi\,D)^{-3/2}\,e^{-\frac{|\x|^2}{4\,D}}  \qquad\mbox{with}\qquad 
          D \geq D_0:= \left(\frac{500}{k_c}\right)^2\,.
$$
Here $D$ is a constant and $\sigma=\sqrt{2\,D}\geq 0.036\cdot mm$ is the variance of $\varphi$. This indicates that 
a resolution of $0.036\cdot mm$ is possible (in the noise-free case).

\section{Conclusions}

In this paper we have proposed a time reversal functional $\I$ for solving PAT of a dissipative medium that obeys the causal wave 
equation of Nachman, Smith and Waag. This model is appropriate for tissue that is similar to water. 
Our theoretical and numerical investigations have shown that 
\begin{itemize}
\item the time reversal image $\I$ for the considered dissipative medium does not give the exact initial pressure function $\varphi$,

\item however, if the compressibility $\kappa_1$ tends to zero, then $\I\to \varphi$ holds. 

\item Moreover, for appropriate conditions, we have 
      $$
          \I\approx 3.5\cdot\varphi
      $$ 
      and a resolution of $\sigma\approx 0.036\cdot mm$ is feasible (in the noise-free case).

\end{itemize}
All in all, it follows that (regularized) time reversal is a valuable solution method for PAT of dissipative 
tissue similar to water.

\section{Appendix}

For the convenience of the reader we list the representations of the solutions $p$ and $q$ of the dissipative 
wave equation~(\ref{weqNSW}) and its time reverse~(\ref{trweqNSW}), respectively, in this appendix. 
We use the following definition of the \emph{Fourier transform}  
$$
     \F\{f\}(\k) := \hat f(\k) := (2\,\pi)^{-3/2}\,\int_{\R^3} e^{\i\,\k\cdot\,\x}\,f(\x)\,\d \x
$$
such that the \emph{convolution theorem} reads as follows
$$
       \F\{f *_\x g\} = (2\,\pi)^{3/2}\,\hat f\,\hat g\,.  
$$
By $k$ we denote the the wave number of the wave vector $\k\in\R^3$, i.e. $k:=|\k|$. 

\subsection*{The solution of the dissipative wave equation}

Fourier transformation of the wave equation~(\ref{weqNSW}) with respect to $\x$ and solving the \emph{Helmholtz equation}  
$$
     -k^2 \left(1 + \tau_1\frac{\partial}{\partial t}\right) \hat p 
     - \frac{1}{c_0^2}\left(1 + \tau_0\frac{\partial}{\partial t}\right)\frac{\partial^2 \hat p}{\partial t^2} 
     = - \frac{\hat\varphi}{c_0^2}\left(1 + \tau_1\frac{\partial }{\partial t}\right)
            \delta'(t) \,,
$$
yields $p = \frac{\partial P}{\partial t}$ with 
$$
   \hat P(\k,t)  
     = \hat \varphi(\k)\,\sum_{j=1}^ 3 A_j\,e^{-\lambda_j(|\k|)\,t}\,H(t) \qquad\mbox{for}\qquad  k\in\R^3 \,,
$$
where $\lambda_j$ and $A_j$ for $j=0,\,1,\,2$ are the solutions of 
\begin{equation}\label{eqlambda}
   -\tau_0\,\lambda^3 + \lambda^2 - c_0^2\,\tau_1\,k^2\,\lambda + c_0^2\,k^2 = 0 
    \qquad\mbox{for}\qquad k\in [0,\infty) \,
\end{equation} 
and
\begin{equation*}\label{eqA}
    \sum_{j=0}^2 A_j\,\lambda_j^m = a_m  \qquad\mbox{for}\qquad m\in\{0,\,1,\,2\}\,
\end{equation*}
with 
\begin{equation*}\label{defam}
        a_0 := 0\,,\qquad
        a_1 := - \frac{\tau_1}{\tau_0}  \qquad\mbox{and}\qquad 
        a_2 := \left(1-\frac{\tau_1}{\tau_0}\right)\, \frac{1}{\tau_0}\,.
\end{equation*} 
The respective \emph{Cardano's formula} read as follows  
\begin{equation*}\label{cardano1}
   \lambda_j =\frac{1}{3\,\tau_0}\,\left( 1 + u_j\,C + \frac{\Delta_0}{u_j\,C}\right) 
             \qquad\mbox{for}\qquad j\in\{0,\,1,\,2\}
\end{equation*}
with 
\begin{equation*}\label{cardano2}
    u_0 = 1\,,\qquad    u_1 = \frac{-1+\i\,\sqrt{3}}{2}\,, \qquad\quad   u_2 = \frac{-1-\i\,\sqrt{3}}{2}\,,
\end{equation*}
\begin{equation*}\label{cardano3}
   \Delta_0 := 1 - 3\,c_0^2\,\tau_0\,\tau_1\,k^2 \,,\qquad\quad 
   \Delta_1 := 2 + 9\,c_0^2\,\tau_0\,(3\,\tau_0 -  \tau_1)\,k^2 \,
\end{equation*}
and
\begin{equation*}\label{cardano4}
         C := \sqrt[3]{\frac{\Delta_1 + \sqrt{\Delta_1^2 - 4\,\Delta_0^3}}{2}}\,.
\end{equation*} 
\begin{figure}[!ht]
\begin{center}
\includegraphics[height=5cm,angle=0]{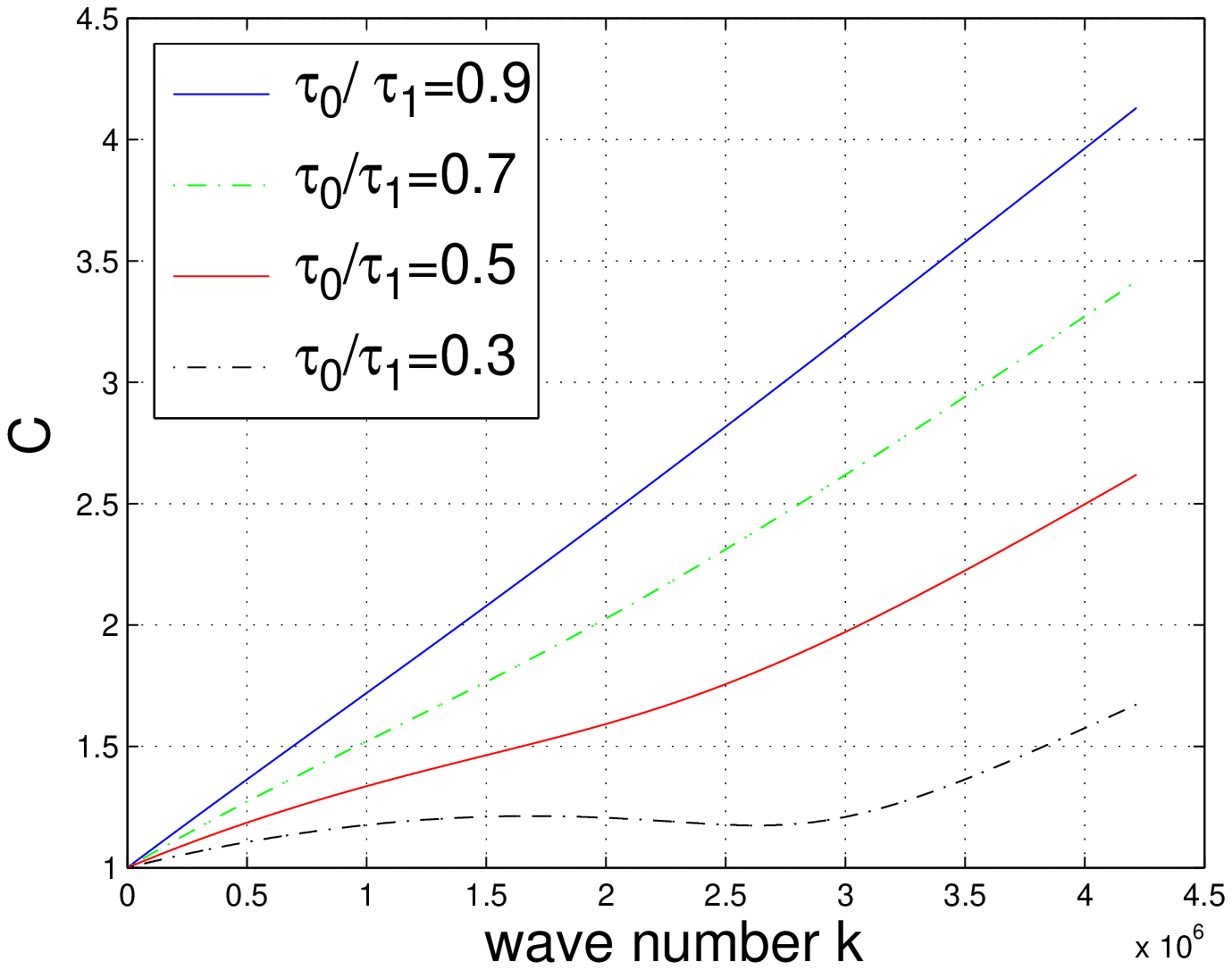}
\includegraphics[height=5cm,angle=0]{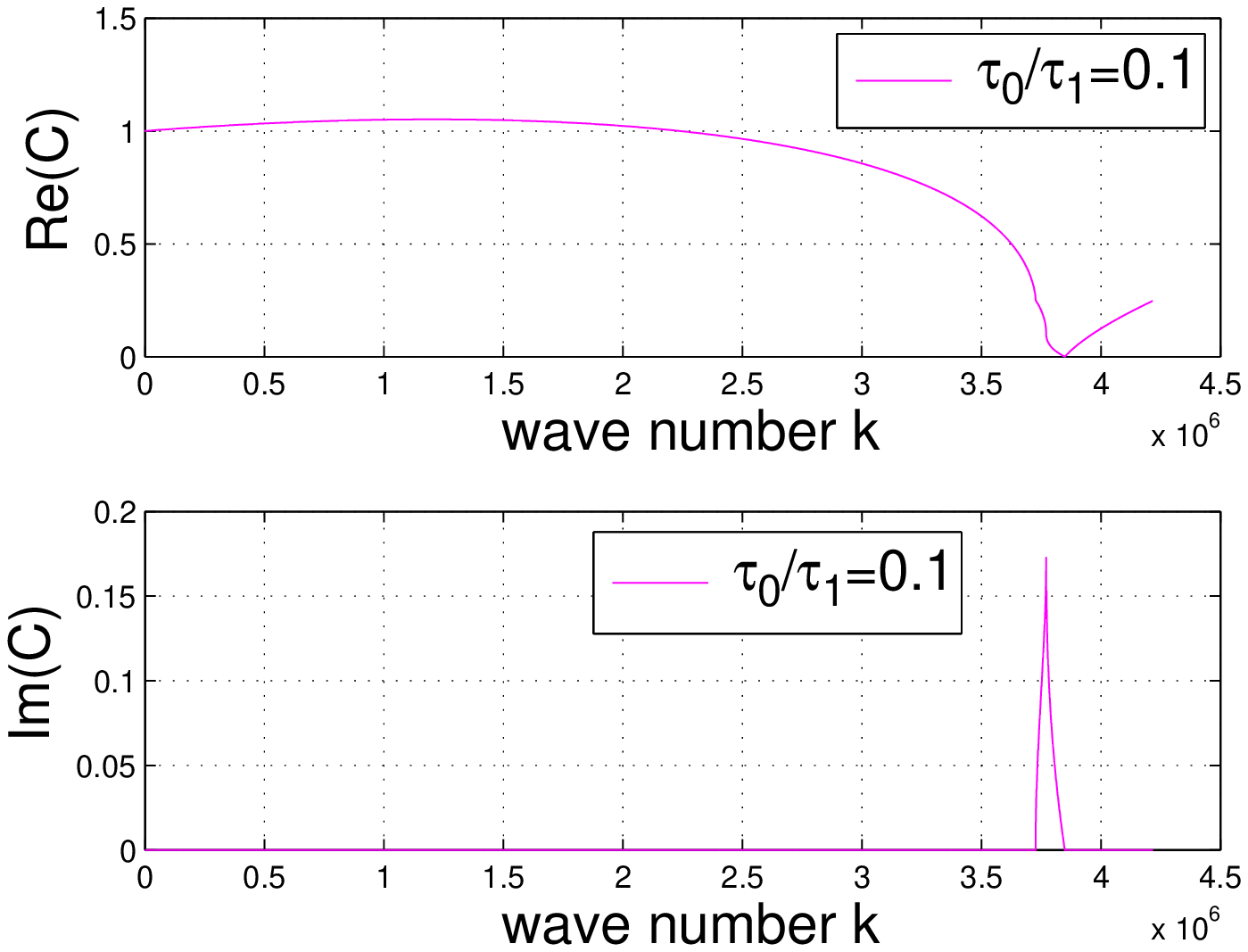}
\end{center}
\caption{Visualization of $C$. We see that $C$ is real-valued for $\frac{\tau_0}{\tau_1}=0.9,\,0.7,\,0.5,\,0.3$ but 
complex-valued for $\frac{\tau_0}{\tau_1}=0.1$. In this paper we focus on tissue similar to water, i.e. 
$\frac{\tau_0}{\tau_1}\approx 0.47$.}
\label{fig:Cs}
\end{figure}
From the above \emph{Cardano's formula}, we get 
\begin{equation*}
\begin{aligned}
   \lambda_0 = \frac{1}{3\,\tau_0}\,\left( 1 + C + \frac{\Delta_0}{C}\right) \,, \qquad \lambda_1= \mu + \i\,\vartheta 
        \qquad\mbox{and}\qquad \lambda_2 = \mu - \i\,\vartheta 
\end{aligned}
\end{equation*}
with
\begin{equation*}
\begin{aligned}
   \mu := \frac{2 - \left(C + \frac{\Delta_0}{C}\right)}{6\,\tau_0} \qquad\mbox{and}\qquad
   \vartheta := \frac{\sqrt{3}\,\left(C - \frac{\Delta_0}{C}\right)}{6\,\tau_0}  \,. 
\end{aligned}
\end{equation*}
For real-valued $C$, it follows that $\lambda_0$, $\mu$ and $\vartheta$ are real-valued. 
For tissue similar to water, we have $\frac{\tau_0}{\tau_1}\approx 0.47$ and thus $C$ is positive 
and real-valued (cf. Fig.~\ref{fig:Cs}). Moreover, it can be shown that
\begin{equation*}
\begin{aligned}
   &A_0 = \frac{a_2-a_1\,(\lambda_2+\lambda_1)}{(\lambda_2-\lambda_0)\,(\lambda_1-\lambda_0)}\,,\qquad
      A_1 = \frac{a_1\,(\lambda_2+\lambda_0) - a_2}{(\lambda_1-\lambda_0)\,(\lambda_2-\lambda_1)}\,\qquad\mbox{and}\\
       &A_2 = \frac{a_2 - a_1\,(\lambda_1+\lambda_0)}{(\lambda_2-\lambda_0)\,(\lambda_2-\lambda_1)} \,.
\end{aligned}
\end{equation*}
If $C$ is real-valued, then it follows from $\lambda_2 = \overline{\lambda_1}$ that $A_0$ is real-valued and 
$A_2 = -\overline{A_1}$ (cf. Fig~\ref{fig:A}), which permits to simplify the representation of the time reversal image 
$\I$ in Section~\ref{sec-opeq}.  
Finally, it can be shown that 
$$
            \lim_{k\to \infty} \left(C+\frac{\Delta_0}{C}\right) = \frac{3\,\tau_0-\tau_1}{\tau_1}
$$ 
and thus
\begin{equation}\label{limmu}
    \lim_{k\to \infty} \lambda_0 = \frac{1}{\tau_1} \qquad\mbox{and}\qquad   
    \lim_{k\to \infty} \mu = \frac{1}{2}\,\left(\frac{1}{\tau_0} - \frac{1}{\tau_1}\right) \,.
\end{equation}
The functions $\lambda_0$, $\mu$, $\vartheta$ and $|\lambda_1|$ for the parameter values as in~(\ref{parvalues}) are 
visualized in Fig.~\ref{fig:01}. For these values, we have 
$\lim_{k\to \infty} \lambda_0\approx 10^9\cdot s^{-1}$ and 
$\lim_{k\to \infty} \mu\approx 5.6250\cdot 10^8\cdot s^{-1}$ (cf. Fig.~\ref{fig:01}).

\subsection*{The solution of the time reversed wave equation}

Similarly as above, it follows for nice $\phi_T$ that the solution of the time reversed wave equation~(\ref{trweqNSW}) 
is given by $q = \frac{\partial Q}{\partial t}$ with
$$
   \hat Q(\k,t)  
      =   - \hat \phi_T(\k)\,\sum_{j=0}^2 A_0\,e^{\lambda_j(|\k|)\,t}\,H(t)\,.
$$
Here $\lambda_j$ and  $A_j$ for $j=0,\,1,\,2$ are defined as above.

\end{document}